\documentclass[procedia]{easychair}

\usepackage{amsmath,amsfonts,amssymb}

\usepackage{doc}
\usepackage{makeidx}

\usepackage{algorithm2e}

\def\procediaConference{99th Conference on Topics of
  Superb Significance (COOL 2014)}

\newcommand{\beq}{\begin{equation}}
\newcommand{\eeq}{\end{equation}}
\newcommand{\mc}{\multicolumn}

\newtheorem{theorem}{\bf Theorem}

\title{A Spectral Projection Preconditioner for Solving Ill Conditioned Linear %
Systems}

\titlerunning{A Spectral Projection Preconditioner}

\author{
    Man-Chung Yeung\inst{1}
\and
    Craig C. Douglas\inst{2,1}
\and
    Long Lee\inst{1}
}

\institute{
University of Wyoming,
  Department of Mathematics,
  Laramie, WY, USA\\
  \email{\{myeung, llee\}@uwyo.edu}
\and
  University of Wyoming,
  School of Energy Resources,
  Laramie, WY, USA\\
  \email{craig.c.douglas@gmail.com}\\
 }

\authorrunning{Yeung, Douglas, and Lee}

\begin{document}

\maketitle


\begin{abstract}
We present a preconditioner based on spectral projection that is combined with a deflated
Krylov subspace method for solving ill conditioned linear systems of
equations.
Our results show that the proposed algorithm requires many fewer iterations to
achieve the convergence criterion for solving an ill conditioned problem than
a Krylov subspace solver.
In our numerical experiments, the solution obtained by the proposed algorithm
is more accurate in terms of the norm of the distance to the exact solution of
the linear system of equations.
\end{abstract}


\section{Introduction}
\label{sect:introduction}

Both the robustness and efficiency of iterative methods are affected by the
condition number of the associated linear system of equations.
When a linear system has a large condition number (usually due to eigenvalues
that are close to the origin of the spectrum domain) iterative methods tend to
take many iterations before a given convergence criterion is satisfied.
Iterative methods may fail to converge within a reasonable elapsed time or
even fail to converge if the condition number is too large.
Unstable linear systems (i.e., systems with large condition numbers) are
called ill conditioned.
For an ill conditioned linear system, slight changes in the coefficient matrix
or the right hand side cause large changes in the solution.
Roundoff errors in the computer arithmetic can cause instability when attempts
are made to solve an ill conditioned system either by direct or iterative
methods on a computer.

It is widely recognized that linear systems resulting from discretizing
ill posed integral equations of the first kind are highly ill conditioned.
This is due to the eigenvalues for the first kind integral equations with
continuous or weakly singular kernels have an accumulation point at zero.
Integral equations of the first kind are frequently seen in statistics, such
as unbiased estimation, estimating a prior distribution on a parameter given
the marginal distribution of the data and the likelihood, and similar tests
for normal theory problems.
They also arise from indirect measurements and nondestructive testing in
inverse problems.
Other ill conditioned linear systems can be seen in training of neural
networks, seismic analysis, Cauchy problem for parabolic equations, and
multiphase flow of chemicals.
For pertinent references of ill conditioned linear systems, see
\cite{engl93, groetsch97}.

Solving ill conditioned linear algebra problems has been a long standing
bottleneck for advancing the use of iterative methods.
The convergence of iterative methods for ill conditioned problems can be
improved by using preconditioning.
Development of preconditioning techniques is therefore a very active research
area.

A preconditioning strategy that {\em deflates} a few isolated external
eigenvalues was first introduced by Nicolaides \cite{nicolaides87} and
investigated by others (e.g., \cite{mansfield91,tv06, vsm99,fv01}).
The deflation strategy is an action that removes the influence of a subspace
of the eigenspace on the iterative process.
A common way to deflate an eigenspace is to construct a proper projector $P$
as a preconditioner and solve
\begin{equation}\label{equ:1-28-1}
PAx= Pb,\quad P, A\in\mathbb{C}^{N\times N}.
\end{equation}
The deflation projector $P$, which is orthogonal to the matrix $A$ and the
vector $b$ against some subspace, is defined by
\begin{equation}\label{eq:projector}
P=I-AZ(Z^{H}AZ)^{-1}Z^{H}, \quad Z\in\mathbb{C}^{N\times m},
\end{equation}
where $Z$ is a matrix of deflation subspace, i.e., the space to be projected
out of the residual and $I$ is the identity matrix of appropriate size
\cite{saad96, fv01}.
We assume that
    (1) $m \ll N$ and
    (2) $Z$ has rank $m$.

A deflated $N{\times}N$ system (\ref{equ:1-28-1}) has an eigensystem different
from that of $Ax =b$.
Suppose that $A$ is diagonalizable.
Set $Z= [v_1,\cdots,v_m]$ whose columns are the eigenvectors of $A$ associated
with the eigenvalues $\lambda_1,\cdots,\lambda_m$.
Then the spectrum
\[
\sigma(PA) = \sigma(A) \setminus \{\lambda_1,\cdots,\lambda_m\}.
\]
The eigenvectors are not easily available, which motivates us to develop an
efficient and robust algorithm for finding an approximate deflation subspace
that does not use the exact eigenvectors to construct the deflation projector
$P$.

Suppose that we want to deflate a set of eigenvalues of $A$ enclosed in a
circle $\Gamma$ that is centered at the origin with the radius $r$.
Without loss of generality, let this set of eigenvalues be
$\{\lambda_1,\cdots,\lambda_k\}$.
Let the subspace spanned by the corresponding eigenvectors of
$\{\lambda_1,\cdots,\lambda_k\}$ be
$\mathcal{Z}_k =\text{Span}\{v_1,\cdots,v_k\}$.
Then the deflation subspace matrix $Z$ in (\ref{eq:projector}) obtained by
randomly selecting $m$ vectors from $\mathcal{Z}_k$ can be written as a
contour integral \cite{saad11}
\beq\label{eq:Z_oint}
Z=\frac{1}{2\pi \sqrt{-1}}\oint_{\Gamma}(zI-A)^{-1}Y\,dz,
\eeq
where $Y$ is a random matrix of size $N{\times}m$.
If the above contour integral is approximated by a Gaussian quadrature, we
have
\beq\label{eq:gaussian}
Z \approx \sum_{i=1}^{q}\omega_i(z_iI-A)^{-1}Y,
\eeq
where $\omega_i$ are the weights, $z_i$ are the Gaussian points, and $q$ is
the number of Gaussian points on $\Gamma$ for the quadrature.

It is worth noting that (\ref{eq:gaussian}) requires the solution of
shifted linear systems $(z_i I-A)X= Y$, $i=1,\cdots,q$.
Each of the $q$ problems can be solved in parallel without communication
before completion.
In addition, each of the $q$ problems can be solved using a parallel solver.
Thus, using (\ref{eq:gaussian}) can lead to the efficient use of many
processors, not just $q$ or $1$.

Using (\ref{eq:gaussian}) for the deflation projector $P$ in
(\ref{eq:projector}), the preconditioned linear system (\ref{equ:1-28-1}) is
no longer severely ill conditioned.
We remark that the construction of a deflation subspace matrix $Z$ through
(\ref{eq:gaussian}) is motivated by the works in \cite{ss03, st07, polizzi,
TP13}.

The outline of the paper is as follows.
In \S\ref{sec:problem}, we introduce some theoretical backgrounds on the
deflated Krylov subspace method (and specifically to GMRES \cite{Saa86S}) and
on the deflation subspace matrix $Z$ in (\ref{eq:Z_oint}).
In \S\ref{Subsec:Problem:NumerExI}, we incorporate the $Z$ computed by
(\ref{eq:gaussian}) into a deflated Krylov subspace method to solve a linear
system.
In \S\ref{Subsec:Problem:NumerExI}, we present numerical experiments.
In \S\ref{sect:future-work}, we briefly introduce state of art parallel
multigrid methods that will be applied to the computation of $Z$ in future.
In \S\ref{sect:conclusions}, we offer conclusions.


\section{Methodology}
\label{sec:problem}

We consider the solution of the linear system
\begin{equation}\label{equ:11-19-1}
A x = b
\end{equation}
by a Krylov subspace method, where we assume that
$A \in {\mathbb C}^{N \times N}$ is nonsingular and $b \in {\mathbb C}^N$.
Let an initial guess $x_0 \in {\mathbb C}^N$ be given and let
$r_0 = b - Ax_0$ be its residual.
A Krylov subspace method recursively constructs an approximate solution,
$x_j$ such that
$$
x_j \in x_0 + {\cal K}_j(A, r_0) \equiv x_0 + \mbox{span} \{r_0, A r_0, \ldots, A^{j-1} r_0\},
$$
where its residual $r_j = b - Ax_j$ satisfies some desired conditions.
The convergence rate of a Krylov subspace method depends on the eigenvalue
distribution of the coefficient matrix $A$.
A variety of error bounds on $r_j$ exist in the literature.
Let us take GMRES \cite{Saa86S} as an example.


\subsection{GMRES}
\label{Subsec:Problem:GMRES}

The residual $r_j$ in the GMRES method is required to satisfy the condition
$$
\|r_j \|_2 = \min_{\xi \in x_0 + {\cal K}_j(A, r_0)}\|b - A \xi\|_2.
$$
Thus, the approximate solution $x_j$ obtained at iteration $j$ of GMRES is
optimal in terms of the residual norm.
In the case where $A$ is diagonalizable, an upper bound on $\| r_j\|_2$ is
provided by the following result.

\begin{figure}[tbh]
\begin{center}
\includegraphics[width=5.8in]{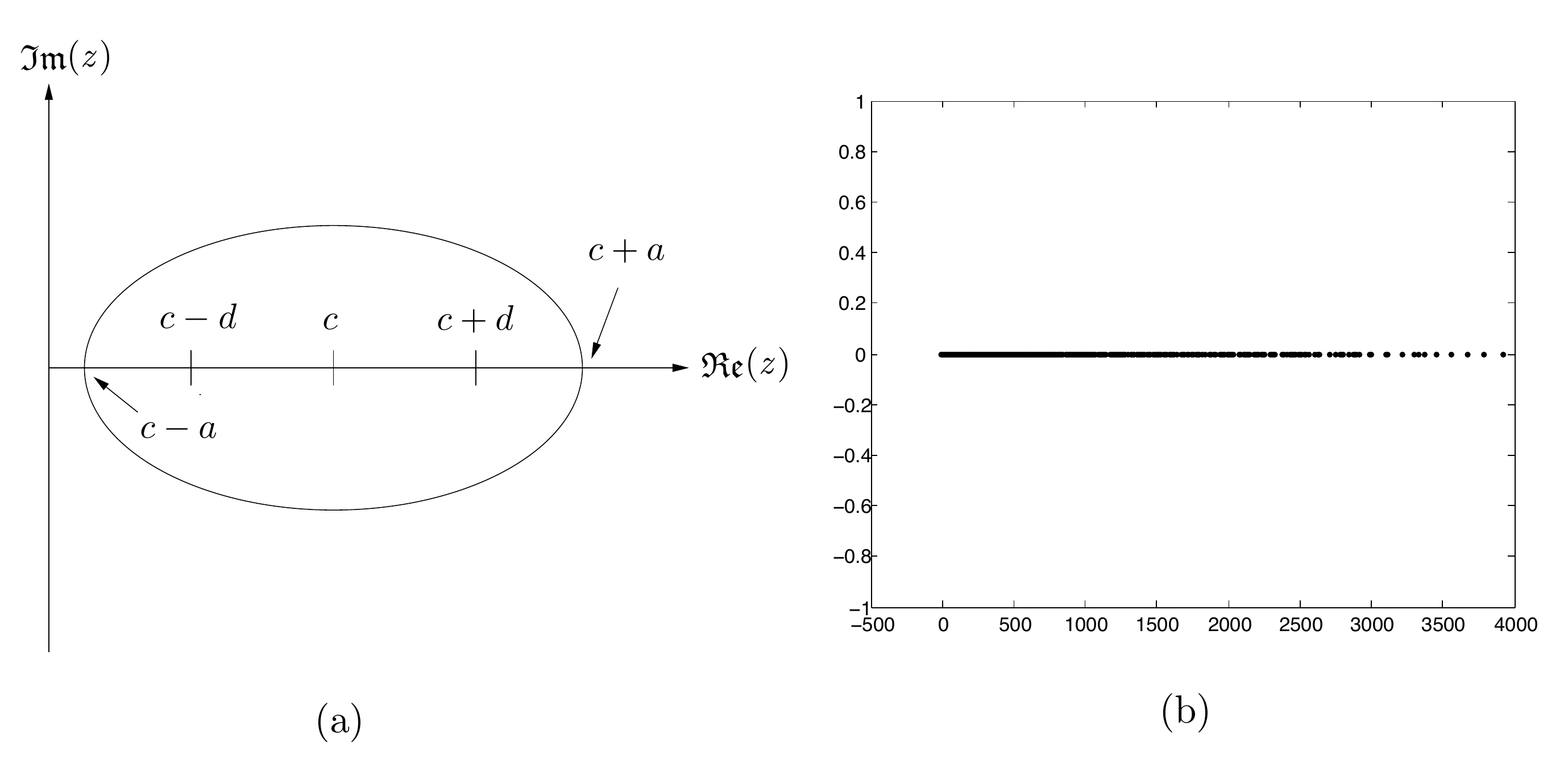}
\end{center}
\vskip -0.25in
\caption{(a) A schematic ellipse in the complex plane with center $c$, focal %
distance $d$, and semi-major axis $a$.%
(b) Eigenvalue distribution of the test matrix {\it bcsstm27}.}
\label{fig:ellipse}
\end{figure}

\begin{theorem}\label{cor:3-16-1} (\cite[Corollary 6.33]{saad96})
Suppose that $A$ can be decomposed as
\begin{equation}\label{equ:11-19-2}
A = V \Lambda V^{-1}
\end{equation}
with $\Lambda$ being the diagonal matrix of eigenvalues.
Let
$E(c, d, a)$ denote the ellipse in the complex plane with center $c$,
focal distance $d$,
and semi-major axis $a$ (see Fig. \ref{fig:ellipse}(a)).
If all the eigenvalues of $A$ are located in $E(c, d, a)$ that excludes the
origin of the complex plane, then
\begin{equation}\label{equ:3-16-31}
\|r_j\|_2 \leq \kappa_2(V) \frac{C_j(\frac{a}{d})}{|C_j(\frac{c}{d})|} \|r_0\|_2,
\end{equation}
where $\kappa_2(V) = \|V\|_2 \|V^{-1}\|_2$ and $C_j$ is the Chebyshev
polynomial of degree $j$.
\end{theorem}

An explicit expression of $C_j(\frac{a}{d})/C_j(\frac{c}{d})$ can be found on
p.~207 of \cite{saad96}, and under some additional assumptions on $E(c, d, a)$
(say, the ellipse in Fig. \ref{fig:ellipse}(a)),
\begin{equation}\label{equ:3-18-153}
\frac{C_j(\frac{a}{d})}{C_j(\frac{c}{d})} \approx \left(\frac{a + \sqrt{a^2 - d^2}}{c + \sqrt{c^2 - d^2}}\right)^j \equiv \delta^j.
\end{equation}

The upper bound in (\ref{equ:3-16-31}) therefore contains two factors: the
condition number $\kappa_2(V)$ of the eigenvector matrix $V$ and the scalar
$\delta$ determined by the distribution of the eigenvalues of $A$.
If $A$ is nearly normal and has a spectrum $\sigma(A)$ which is clustered
around $1$, we would have $\kappa_2(V) \approx 1$ and $\delta < 1$.
In this case, $\|r_j\|_2$ decays exponentially in a rate of power $\delta^j$,
resulting in a fast convergence of GMRES.
The error bound (\ref{equ:3-16-31}) hides the fact that the convergence rate
is better if the eigenvalues of $A$ are clustered \cite{vdS86v}.


Since the ellipse $E(c, d, a)$ in Theorem~\ref{cor:3-16-1} must include all of
the eigenvalues of $A$, the outlying eigenvalues may keep the ellipse large,
implying a large $\delta$.
To reduce $\delta$, we therefore only remove these outlying eigenvalues from
$\sigma(A)$.
Any procedure of doing so is known as {\it deflation}.
GMRES in combination with deflation is called {\it Deflated GMRES}.


\subsection{Deflated GMRES}
\label{sec:conv}

Suppose $x^*$ is the exact solution of (\ref{equ:11-19-1}).
Let a so-called deflation-subspace matrix
$Z = [z_1, \ldots, z_m] \in {\mathbb C}^{N \times m}$ be given, whose columns
are linearly independent.
Define the two projectors \cite{dinkla, fv01, vsm99, YTV}
\begin{equation} \label{equ:8-21-1}
\begin{array}{rcl} P \equiv I - A Z M^{-1} Z^H &
\mbox{and} & \widetilde{P} \equiv I - Z M^{-1} Z^H A,
\end{array}
\end{equation}
where $M = Z^H A Z$ is assumed to be invertible. It is straightforward to
verify that $P^2 = P, \; \widetilde{P}^2 = \widetilde{P}$ and
$P A = A \widetilde{P}$.

Using $\widetilde{P}$, we split $x^*$ into two parts:
$$
x^* = (I - \widetilde{P}) x^* + \widetilde{P} x^* \equiv x^*_1 + x^*_2.
$$
For $x^*_1$, we have
$$
x^*_1 = (I - \widetilde{P}) x^* = Z M^{-1} Z^H A x^* = Z M^{-1} Z^H b.
$$
For
$x^*_2$, we obtain
$$
x^*_2 = A^{-1} P b,
$$
since
$
A x^*_2 = A \widetilde{P} x^* = P A x^* = P b.
$
Now, if $x^\#$ is a solution of the singular system
\begin{equation}\label{equ:7-31}
P A x = P b,
\end{equation}
then
$$\begin{array}{rcl}
A \widetilde{P} x^\# = P b & \Leftrightarrow & \widetilde{P} x^\# = A^{-1} P b
 = x^*_2.
\end{array}
$$
Based on the above observation, a Deflated GMRES algorithm for the solution of
(\ref{equ:11-19-1}) is given in Algorithm~\ref{alg:DGMRES}.

\begin{algorithm}
  \SetKwData{Left}{left}\SetKwData{This}{this}\SetKwData{Up}{up}
  \SetKwFunction{Union}{Union}\SetKwFunction{FindCompress}{FindCompress}
  \SetKwInOut{Input}{input}\SetKwInOut{Output}{output}
\begin{tabbing}
x\=xxx\= xxx\=xxx\=xxx\=xxx\=xxx\kill
\>1.\>Choose $Z$; \\
\>2.\>Compute $x_1 = Z M^{-1} Z^H b$; \\
\>3.\>Solve $P A x = P b$ by GMRES to obtain a solution $x^\#$; \\
\>4.\>Compute $x_2 = \widetilde{P} x^\#$; \\
\>5.\>Compute $x = x_1 + x_2$.
\end{tabbing}
\caption{Deflated GMRES}
\label{alg:DGMRES}
\end{algorithm}

We remark that the GMRES in Algorithm \ref{alg:DGMRES} can be replaced by any
other linear solvers, direct or iterative.
We note that when $A$ is symmetric, indefinite that we can use the MINRES
method \cite{Minres} instead of GMRES and compute the same solution $x$ using
far less memory.

Assume that the nonsingular $A \in {\mathbb C}^{N \times N}$ has a
decomposition~(\ref{equ:11-19-2}) with $V = [v_1, \ldots, v_N]$ and
$\Lambda = \mbox{diag}\{\lambda_1, \ldots, \lambda_N\}$.
If we set $Z = [v_1, \ldots, v_m]$ in (\ref{equ:8-21-1}), then the spectrum
$\sigma(PA)$ contains all the eigenvalues of $A$ except
$\lambda_1, \ldots, \lambda_m$, namely,
$\sigma(PA) = \{0, \cdots, 0, \lambda_{m+1}, \cdots, \lambda_N\}$.

Perform a $QR$ factorization on $V$ as follows:
$$
V = Q R \equiv [Q_1, Q_2] \left[
\begin{array}{cc}
R_{11} & R_{12}\\
0 & R_{22}
\end{array} \right],
$$
where $Q_1 \in {\mathbb C}^{N \times m}$ and
$R_{11} \in {\mathbb C}^{m \times m}$.
If we set $Z = [v_1, \ldots, v_m]$ and apply GMRES to solve (\ref{equ:7-31}),
an upper bound on $\|r_j\|$ is given by the following theorem \cite{YTV}.

\begin{theorem}\label{cor:4-4-50} Suppose that $A$ has a decomposition
(\ref{equ:11-19-2}) and GMRES is used to solve (\ref{equ:7-31}) with
$Z = [v_1, \ldots, v_m]$.  
If all the eigenvalues $\lambda_{m+1}, \ldots, \lambda_N$ of $A$ are located
in an ellipse $E(c, d, a)$ that excludes the origin of the complex plane, then
\begin{equation}\label{equ:3-18-32}
\|r_j\|_2 \leq \kappa_2(R_{22}) \frac{C_j(\frac{a}{d})}{|C_j(\frac{c}{d})|} \|r_0\|_2. 
\end{equation}
\end{theorem}

With (\ref{equ:3-18-153}), the upper bound (\ref{equ:3-18-32}) of the residual
norm $\|r_j\|_2$ of
GMRES
is determined by the condition number of $R_{22}$ (rather than $V$) and the
scalar $\delta$ which is determined by the distribution of the undeflated
eigenvalues $\lambda_{m+1}, \ldots, \lambda_N$ of $A$.
The $\delta$ in Theorem~\ref{cor:4-4-50} is generally less than the $\delta$
associated with Theorem~\ref{cor:3-16-1}.
This partially explains why an eigenvalue-deflation is likely to lead to a
faster convergence of GMRES.

The proof of Theorem \ref{cor:4-4-50} is based on the observation that when
GMRES solves a singular linear system it is actually solving a nonsingular
linear system of smaller size.
Theorem \ref{cor:4-4-50} then follows from the application of Theorem
\ref{cor:3-16-1} to the nonsingular linear system.
See \cite{YTV} for the details of the proof of Theorem \ref{cor:4-4-50}.


\subsection{Spectral Projector and Construction of $Z$}
\label{Subsec:Problem:SpectralProj}

A spectral projector is described in detail in \S3.1.3-\S3.1.4 of
\cite{saad11}.
Other references include \cite{chatelin, dunford, kato}.
Let $A = VJV^{-1}$ be the Jordan canonical decomposition of $A$, where
$$\begin{array}{ccc}
V = [v_1, v_2, \ldots, v_N] & \mbox{and} & J = diag\{J_{N_1}(\lambda_1), J_{N_2}(\lambda_2), \ldots, J_{N_d}(\lambda_d)\}.
\end{array}
$$
The diagonal block $J_{N_i}(\lambda_i)$ in $J$ is an $N_i \times N_i$ Jordan
block associated with the eigenvalue $\lambda_i$.
The eigenvalues $\lambda_i$ are not necessarily distinct and can be repeated
according to their multiplicities.

Let $\Gamma$ be a given positively oriented simple closed curve in the complex
plane.
Without loss of generality, let the set of eigenvalues of $A$ enclosed by
$\Gamma$ be $\{\lambda_1, \lambda_2, \ldots, \lambda_k\}$ so that the
eigenvalues $\lambda_{k+1}, \ldots, \lambda_d$ lie outside the region enclosed
by $\Gamma$.
Set $s \equiv N_1 + N_2 + \ldots + N_k$, the number of eigenvalues inside
$\Gamma$ with multiplicity taken into account.
Then the residue
$$
P_\Gamma = \frac{1}{2 \pi \sqrt{-1}} \oint_\Gamma (zI - A)^{-1} dz
$$
defines a projection operator onto the space
$\sum_{i = 1}^k \mbox{Null}(A - \lambda_i I)^{l_i}$, where $l_i$ is the index
of $\lambda_i$, namely,
$$
\mbox{Range}(P_\Gamma) = \mbox{span}\{v_1, v_2, \ldots, v_s\}.
$$
In particular, if $A$ has a diagonal decomposition (\ref{equ:11-19-2}),
$P_\Gamma$ is a projector onto the sum
$\sum_{i = 1}^k {\mathbb E}_{\lambda_i}$ of the $\lambda_i$-eigenspace
${\mathbb E}_{\lambda_i}$ of $A$.

Pick a random matrix $Y \in {\mathbb C}^{N \times m}$ and set
\begin{equation}\label{equ:11-24-1}
Z = P_\Gamma Y = \frac{1}{2 \pi \sqrt{-1}} \oint_\Gamma (zI - A)^{-1} Y dz
\end{equation}
in (\ref{equ:8-21-1}).
In the case where $m = s$, we almost surely have $\sigma(PA) = \{0, \cdots, 0,
\lambda_{k+1}, \cdots, \lambda_d\}$.
Therefore all the eigenvalues of $A$ inside $\Gamma$ are removed from the
spectrum of $PA$.


\section{Numerical Examples}
\label{Subsec:Problem:NumerExI}

In this section, we demonstrate the effect of the deflation-subspace matrix
$Z$ defined by (\ref{equ:11-24-1}) and computed by the Legendre-Gauss
quadrature on the solution of the linear system (\ref{equ:11-19-1}).

All the computations were done in Matlab Version 7.1 on a Windows 10 machine.
Besides GMRES, we also employ a modified BiCG (MBiCG) as the Krylov solver in
Line 3 of Algorithm~\ref{alg:DGMRES} and for finding the solution of all of
the linear systems in the computation of $Z$.
The MBiCG is the standard BiCG, but outputs the computed solution $x$ that
either satisfies {\it relres} $<$ {\it tol} or has the smallest relative
residue among all the computed $x$ from iteration $1$ to iteration
{\it maxit}, where {\it relres} is the relative residue of $x$, {\it tol} is
the user supplied input stopping tolerance, and {\it maxit} is the maximum
number of iterations.

The
numerical solution using a deflated, restarted GMRES of the linear systems
obtained from the discretization of the two dimensional steady-state
convection-diffusion equation
\begin{equation}\label{equ:11-24-10}
\begin{array}{rcl}
- u_{xx} - u_{yy} - Re\, (p(x, y) u_x - q(x, y) u_y)] = f(x, y), & & (x, y) \in [0, 1]^2
\end{array}
\end{equation}
with Dirichlet boundary conditions was studied in depth in \cite{dinkla}.
This steady-state version (\ref{equ:11-24-10}) of the two-dimensional
convection-diffusion equation is from \cite{zhang}.

In \cite{dinkla}, two types of deflation-subspace matrices $Z$ were used:
eigenvectors obtained from the Matlab function $eig$ and algebraic subdomain
vectors.
The $Z$ of algebraic subdomain vectors works well for the fluid flow problem
(\ref{equ:11-24-10}), but seems not for the other problem presented in \cite{dinkla}.
Accurately calculating eigenvalues of large matrices is very time consuming.
Therefore deflation with the $Z$ of true eigenvectors is not practicable.

Numerical experiments in \cite{dinkla} have shown that eigenvalues close to
the origin hamper the convergence of a Krylov subspace method.
Hence, deflation of these eigenvalues is very beneficial.
Based on this observation, we chose in our experiments the integration path
$\Gamma$ in (\ref{equ:11-24-1}) to be a circle $D(c, r)$ with the center $c$
near the origin.
For the $Y$ in (\ref{equ:11-24-1}), we picked a random $Y \in {\mathbb R}^{N
\times m}$ by the Matlab command $Y = randn(N,m)$ with $m$ not less than the
exact number $s$ of eigenvalues inside $\Gamma$.
We remark that an efficient stochastic estimation method of $s$ has been
developed in \cite{futa}.

We computed the integral in (\ref{equ:11-24-1}) by the
Legendre-Gauss quadrature
\begin{equation}\label{equ:11-24-13}
Z = \frac{r}{2} \int_{-1}^1 e^{\pi \theta i} ((c+r e^{\pi \theta i}) I - A)^{-1} Y d \theta
 \approx \frac{r}{2} \sum_{k = 1}^q \omega_k e^{\pi \theta_k i} ((c+r e^{\pi \theta_k i}) I - A)^{-1} Y,
\end{equation}
where $i = \sqrt{-1}$, and $\omega_k$ and $\theta_k$ are the Legendre-Gauss
weights and nodes on the interval $[-1, 1]$ with truncation order $q$.

In (\ref{equ:11-24-13}), there are $mq$ linear systems
\begin{equation}\label{equ:8-17-1}
\begin{array}{rccl}
((c+r e^{\pi \theta_k \sqrt{-1}}) I - A) x = y_j, & & k = 1, \ldots, q, & j = 1, \ldots, m
\end{array}
\end{equation}
to solve.
We solved all of them by GMRES or MBiCG with the stopping tolerance
$tol = 10^{-15}$ and the maximum number of iterations $maxit = 500$ or $1000$.

Mathematically, the rank of $Z$ defined by (\ref{equ:11-24-1}) is less than
$m$ when $m > s$.
As a result, the matrix $M$ in (\ref{equ:8-21-1}) is singular.
In practice, $Z$ is approximated by (\ref{equ:11-24-13}) and the matrix $M$
becomes near-singular.
In order to remedy this difficulty, one can use the singular value
decomposition (SVD) or the QR decomposition of $Z$ to detect and remove its
nearly dependent columns.
The SVD or QR decomposition involves a high communication cost and may not be
favorable for a parallel computation.
In this paper, we suggest a column deflation mechanism based on Gaussian
elimination (GE) with complete pivoting \cite{gvl} to remove the dependent
columns in $Z$.

The rationale of the mechanism is as follows.
Let $\widehat{Z} = Z^H Z \in {\mathbb C}^{m \times m}$.
It can be seen that $rank(Z) = rank(\widehat{Z})$.
We perform GE with complete pivoting on $\widehat{Z}$ to reduce $\widehat{Z}$
into an upper triangular form.
Accordingly, we interchange the corresponding columns in $Z$.
If at some step of $j$ the lower right block $\widehat{Z}_{(j+1:m,j+1:m)} =
0$, then the rank of $\widehat{Z}$ is $j$ and $Z_{(:,1:j)}$ consists of all
the linearly independent columns of $Z$.
The purpose that we left-multiply $Z$ by $Z^H$ to form $\widehat{Z}$ is
(i) to reduce the row size of $Z$ from $N$ to $m$, and
(ii) to reduce data movements in the GE process.
See Algorithm~\ref{alg:8-18-1} in \S\ref{sect:appendix} for
implementation details.

In our experiments, we performed the following eight computations.
Numerical results are summarized in Tables \ref{tab:9-3-1}-\ref{tab:8-24-2}.
\begin{enumerate}
\item[\#1.]
Solve (\ref{equ:11-19-1}) by GMRES (or MBiCG) with the initial guess
$x = 0$ with $tol = 10^{-7}$ and $maxit = 10^3 N$.
Compute the true relative errors {\it relres2} $=\|b - A x\|_2/\|b\|_2$ and
{\it relerr} $=\|x - x^*\|_2 / \|x^*\|_2$, where $x$ is the computed solution
output by GMRES (or MBiCG) and $x^*$ is the exact solution of
(\ref{equ:11-19-1}).
\item[\#2.]
Compute by the Matlab function {\it eig} the eigenvectors $v_1, \ldots, v_s$
of $A$ whose associated eigenvalues lie inside $\Gamma$.
Set $Z = [v_1, \ldots, v_s]$.
Perform Algorithm \ref{alg:DGMRES} with GMRES (or MBiCG) as its linear solver
(see Line 3 of Algorithm \ref{alg:DGMRES}).
Set the initial guess $x = 0$ with $tol = 10^{-7}$ and $maxit = 10N$ for GMRES
(or MBiCG).
Compute the true relative errors
{\it relres1} $=\|Pb - PAx^\#\|_2/\|Pb\|_2$,
{\it relres2} $=\|b - A x\|_2/\|b\|_2$, and
{\it relerr} $=\|x - x^*\|_2 / \|x^*\|_2$, where $x^\#$ and $x$ are the
computed solutions in Lines 3 and 5 of Algorithm \ref{alg:DGMRES},
respectively.
\item[\#3.]
Solve all the linear systems in (\ref{equ:8-17-1}) by GMRES (or MBiCG) with
the initial guess $x = 0$ with $tol = 10^{-15}$ and $maxit = 500$.
Compute $Z$ using (\ref{equ:11-24-13}).
Perform Algorithm \ref{alg:DGMRES}
with the initial guess $x = 0$ with $tol = 10^{-7}$ and $maxit = 10N$ for
GMRES (or MBiCG).
Compute the true relative errors {\it relres1}, {\it relres2}, and
{\it relerr} defined in item \#2.
\item[\#4.]
Perform the same computations as described in item \#3 with
$maxit = 1000$ instead of $maxit = 500$.
\item[\#5.]
Solve all the linear systems in (\ref{equ:8-17-1}) by GMRES (or MBiCG) with
the initial guess $x = 0$ with $tol = 10^{-15}$ and $maxit = 500$.
Compute $Z$ using (\ref{equ:11-24-13}).
Setting $\alpha = 10^{-8}$ and $tol\_cge = 10^{-2}$, perform Algorithm
\ref{alg:8-18-1} on the computed $Z$ to remove its nearly dependent columns.
Using the $Z$ output from Algorithm \ref{alg:8-18-1}, perform Algorithm
\ref{alg:DGMRES}
with the initial guess $x = 0$ with $tol = 10^{-7}$ and $maxit = 10N$ for
GMRES (or MBiCG).
Compute the true relative errors {\it relres1}, {\it relres2}, and
{\it relerr} defined in item \#2.
\item[\#6.]
Perform the same computations as described in item \#5 with
$maxit = 1000$ instead of $maxit = 500$.
\item[\#7.]
Solve all the linear systems in (\ref{equ:8-17-1}) by GMRES (or MBiCG) with
the random initial guess $x = randn(N,1)$ with $tol = 10^{-15}$, and
$maxit = 500$.
Compute $Z$ through (\ref{equ:11-24-13}).
Perform Algorithm \ref{alg:DGMRES} with the initial guess $x = 0$ with $tol =
10^{-7}$, and $maxit = 10N$ for GMRES (or MBiCG).
Compute the true relative errors {\it relres1}, {\it relres2}, and
{\it relerr} defined in item \#2.
\item[\#8.]
Perform the same computations as described in item \#7 with
$maxit = 1000$ instead of $maxit = 500$.
\end{enumerate}


\subsection{Example 1}
\label{Subsec:NumExp:Example1}

As in \cite{dinkla}, consider the convection-diffusion equation
(\ref{equ:11-24-10}) with Dirichlet boundary conditions.
The convection coefficients $p(x, y)$, $q(x, y)$, and the source term
$f(x, y)$ were chosen as
$$
\begin{array}{l}
p(x, y) = - \sin x \cos \pi y\\
q(x, y) = \cos(\pi x) \sin y,\\
f(x,y) = 52 \sin (4x+6y) - (4p(x,y)-6q(x,y)) \cos(4x+6y).
\end{array}
$$
Equation (\ref{equ:11-24-10}) was discretized on the unit square $[0,1]^2$
using a $5$-point central difference scheme with a uniform mesh size of
$h = 1/100$ with resulting linear systems (\ref{equ:11-19-1}) of size
$N = 99^2$.

The Reynolds number $Re$ controls the degree of the nonsymmetry in the
coefficient matrix $A$ of (\ref{equ:11-19-1}).
When $Re = 0$, $A$ is symmetric.
As $Re$ increases the nonsymmetry in $A$ also increases.

\begin{figure}[tbh]
\begin{center}
\includegraphics[width=5.8in]{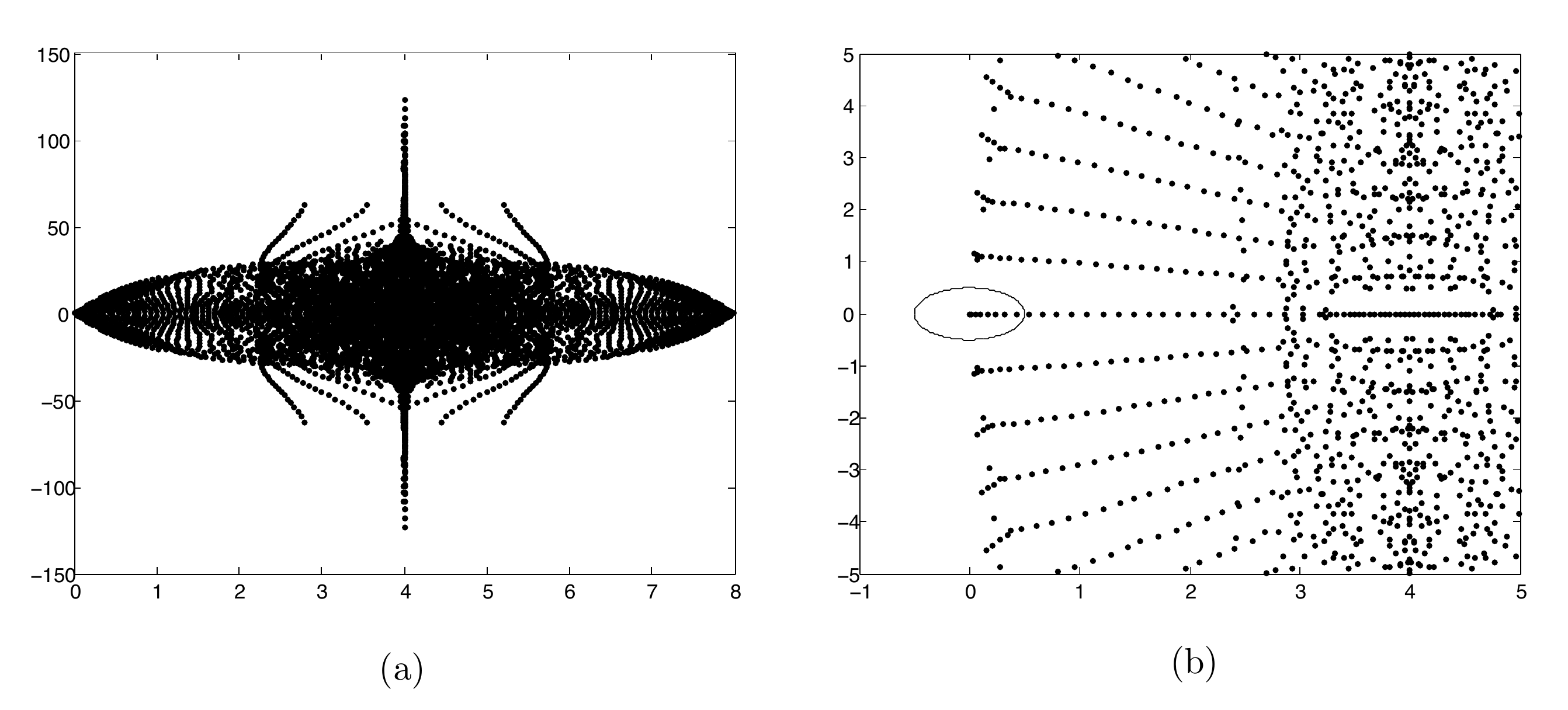}
\end{center}
\vskip -0.25in
\caption{(Example 1)
(a) Eigenvalue distribution of the test matrix $A$. %
(b) Eight eigenvalues of $A$ lie inside $\Gamma$, the circle centered at %
the origin with radius $0.5$. %
The eigenvalue with the smallest distance to the origin is $4.3e\!-\!3$.}
\label{fig:conv_diff}
\end{figure}

In our experiments, we pick the coefficient matrix $A$, but set the
right hand side $b = A {\bf 1}$, where ${\bf 1} = [1, 1, \ldots, 1]^T$.
We know {\em a priori} the exact solution of (\ref{equ:11-19-1}) and the
relative error {\it relerr} is computable.
We choose $Re = 8000$.
The integration path $\Gamma$ is the circle whose center is the origin and
has radius $0.5$.
Figure~\ref{fig:conv_diff} presents the eigenvalue distribution of $A$.

Comparing the numerical results of Computations \#1 and \#2 in Table
\ref{tab:9-3-1}, we see that the convergence of GMRES can be made much faster
with an appropriate eigenvalue-deflation.
Specifically GMRES takes $3295$ iterations to solve (\ref{equ:11-19-1}), but
only $1815$ iterations to solve (\ref{equ:7-31}).
This experiment supports the observation made in \cite{dinkla} that
eigenvalues close to the origin hamper the convergence of a Krylov subspace
method.
It also supports the theory presented in \S\ref{sec:conv}.

The behavior of a Krylov subspace method on the solution of (\ref{equ:7-31})
depends on
(i) the quality of the approximation of the computed $Z$ to the exact
eigenvectors,
(ii) the column size $m$ of $Z$,
(iii) the condition of $M$, and
(iv) perhaps some other unknown factors.

Consider the computed $Z$'s in Tables \ref{tab:9-3-2} and \ref{tab:9-3-3}.
Comparing the corresponding relative residue ranges in the second columns of
the tables we see that the $Z$'s in Table~\ref{tab:9-3-3} are closer to
the exact eigenvectors than the $Z$'s in Table~\ref{tab:9-3-2} are since the
range intervals in Table~\ref{tab:9-3-3} are closer to the origin $0$.
As a result, the GMRES and MBiCG results in Table~\ref{tab:9-3-3} converge
faster than those in Table~\ref{tab:9-3-2} when solving (\ref{equ:7-31}).

Consider Table~\ref{tab:9-3-3}.
As $m$ is increases from $10$ to $50$, the condition numbers Cd($M$) of $M$
does not increase substantially.
In this case, the number of iterations needed by GMRES and MBiCG to converge
decreases (see the 5th column of the table).

In Example~1, the exact number $s$ of eigenvalues of $A$ inside the circle
$\Gamma$ is $8$.
The rank of $Z$ in (\ref{equ:11-24-1}) is $8$ mathematically and therefore
the $2$-norm condition number of $Z$ is $\infty$ in the case when $m > s$.
However, the condition numbers of the computed $Z$'s in Tables \ref{tab:9-3-2}
and \ref{tab:9-3-3} are small finite numbers.
This implies that those computed $Z$'s are inaccurate.
They only contain partial information about the true eigenvectors.
Impressively, the computed $Z$'s in Table~\ref{tab:9-3-3} associated with
$m = 50$ perform even better than the true eigenvectors do (compare the 4th
column in Table~\ref{tab:9-3-1} and the 5th column in Table~\ref{tab:9-3-3}).

Since the computed $Z$'s in Example~1 perform so well, there is no reason to
apply Algorithm~\ref{alg:8-18-1} to improve their conditions.

\begin{table}
\centering
\footnotesize{
{\newcommand{\q}[1]{\mc{1}{|l|}{\small\tt #1}}
\noindent
\begin{tabular}{*{7}{|c}|} \hline
\mc{3}{|c|}{Computation \#1}& \mc{4}{|c|}{Computation \#2}
\\ \hline\hline
 \mc{3}{|c|}{ Solve (\ref{equ:11-19-1}) by GMRES} &\mc{2}{|c|}{Solve (\ref{equ:7-31}) by GMRES} & \mc{2}{|c|}{Algorithm \ref{alg:DGMRES}} \\ \hline
 \#iter & relres2 & relerr & \#iter & relres1 & relres2 & relerr
\\ \hline
 $3295$ &$9.9e\!-\!8$ &$3.3e\!-\!7$ &$1815$  &$9.9e\!-\!8$ &$9.9e\!-\!8$ &$3.4e\!-\!6$
\\ \hline\hline
\mc{3}{|c|}{ Solve (\ref{equ:11-19-1}) by MBiCG} & \mc{2}{|c|}{Solve (\ref{equ:7-31}) by MBiCG} &\mc{2}{|c|}{Algorithm \ref{alg:DGMRES}}
\\ \hline
\#iter & relres2 & relerr & \#iter & relres1 & relres2 & relerr
\\ \hline
$10125$ &$9.5e\!-\!8$  &$2.7e\!-\!8$  & $3951$ &$9.9e\!-\!8$ &$9.9e\!-\!8$ &$2.2e\!-\!8$
\\ \hline
\end{tabular}}
}
\vskip -0.10in
\caption{(Example 1) $\Gamma$ is the circle centered at the origin with radius
$0.5$. %
The number of eigenvalues of $A$ inside $\Gamma$ is $8$. %
{\it \#iter} stands for the number of iterations needed to converge.}
\label{tab:9-3-1}
\end{table}

\begin{table}
\centering
\footnotesize{
{\newcommand{\q}[1]{\mc{1}{|l|}{\small\tt #1}}
\noindent
\begin{tabular}{*{8}{|c}|} \hline
\mc{8}{|c|}{Computation \#3}
\\ \hline\hline
\mc{2}{|c|}{Solve (\ref{equ:8-17-1}) by GMRES} & \mc{2}{|c|}{Compute $Z$ by (\ref{equ:11-24-13})} &\mc{2}{|c|}{Solve (\ref{equ:7-31}) by GMRES} & \mc{2}{|c|}{Algorithm \ref{alg:DGMRES}}
\\ \hline
$m$& Rel. Res. Range & Cd($Z$) & Cd($M$)& \#iter & relres1 & relres2 & relerr
\\ \hline 
 $10$&$[1.4e\!-\!3, 1.1e\!-\!1]$  &$6.7$ &$1.8e\!+\!1$ &$2616$  &$9.9e\!-\!8$ &$1.0e\!-\!7$ &$1.9e\!-\!6$ \\ \hline
 $50$ & $[1.1e\!-\!3, 1.2e\!-\!1]$ &$4.5e\!+\!1$ &$3.0e\!+\!2$  &
 $1321$
 &$9.9e\!-\!8$ &$9.9e\!-\!8$  &$9.2e\!-\!7$
\\ \hline\hline
\mc{2}{|c|}{Solve (\ref{equ:8-17-1}) by MBiCG} & \mc{2}{|c|}{Compute $Z$ by (\ref{equ:11-24-13})} &\mc{2}{|c|}{Solve (\ref{equ:7-31}) by MBiCG} & \mc{2}{|c|}{Algorithm \ref{alg:DGMRES}}
\\ \hline
$m$& Rel. Res. Range & Cd($Z$) & Cd($M$)& \#iter & relres1 & relres2 & relerr
\\ \hline 
 $10$& $[2.6e\!-\!2, 5.5e\!-\!1]$ &$5.9e\!+\!2$ & $1.8e\!+\!4$& $6875$ & $6.8e\!-\!8$ & $8.4e\!-\!8$ & $6.2e\!-\!8$\\ \hline
 $50$ & $[2.3e\!-\!2, 5.2e\!-\!1]$ & $3.0e\!+\!6$& $6.7e\!+\!10$ &$6078$ & $9.9e\!-\!8$& $1.4e\!-\!7$ & $3.1e\!-\!7$
\\ \hline
\end{tabular}}
}
\vskip -0.10in
\caption{(Example 1) The integration path $\Gamma$ is the same as in %
Table~\ref{tab:9-3-1}. %
The $q$ in (\ref{equ:11-24-13}) was set to be $q = 2^4$. %
After each of the linear systems in (\ref{equ:8-17-1}) was solved, %
the true relative residue %
$\|y_j - ((c + r e^{\pi \theta_k \sqrt{-1}})I-A)x\|_2/\|y_j\|_2$ %
of the approximate solution $x$ was computed. %
The intervals in the column titled ``Rel. Res. Range'' are the smallest %
intervals that contain these relative residues. %
{\it Cd$(Z)$} and {\it Cd$(M)$} stand for the $2$-norm condition numbers of %
$Z$ and $M$ respectively.}
\label{tab:9-3-2}
\end{table}

\begin{table}
\centering
\footnotesize{
{\newcommand{\q}[1]{\mc{1}{|l|}{\small\tt #1}}
\noindent
\begin{tabular}{*{8}{|c}|} \hline
\mc{8}{|c|}{Computation \#4}
\\ \hline\hline
 \mc{2}{|c|}{Solve (\ref{equ:8-17-1}) by GMRES} & \mc{2}{|c|}{Compute $Z$ by (\ref{equ:11-24-13})} &\mc{2}{|c|}{Solve (\ref{equ:7-31}) by GMRES} & \mc{2}{|c|}{Algorithm \ref{alg:DGMRES}}
\\ \hline
$m$& Rel. Res. Range & Cd($Z$) & Cd($M$)& \#iter & relres1 & relres2 & relerr
\\ \hline 
 $10$ &$[4.1e\!-\!6, 6.0e\!-\!2]$ &$1.1e\!+\!1$ &$8.0e\!+\!1$&$2420$ &$9.9e\!-\!8$ &
$9.9e\!-\!8$ &$2.4e\!-\!6$
\\ \hline
 $50$ & $[3.0e\!-\!6, 8.5e\!-\!2]$ &$1.3e\!+\!2$  &$1.9e\!+\!3$ &$1340$  &$9.9e\!-\!8$  & $9.9e\!-\!8$ &$1.8e\!-\!6$
\\ \hline\hline
 \mc{2}{|c|}{Solve (\ref{equ:8-17-1}) by MBiCG} & \mc{2}{|c|}{Compute $Z$ by (\ref{equ:11-24-13})} &\mc{2}{|c|}{Solve (\ref{equ:7-31}) by MBiCG} & \mc{2}{|c|}{Algorithm \ref{alg:DGMRES}}
\\ \hline
$m$& Rel. Res. Range & Cd($Z$) & Cd($M$)& \#iter & relres1 & relres2 & relerr
\\ \hline 
 $10$ & $[9.9e\!-\!4, 5.8e\!-\!1]$ & $9.6$& $3.3e\!+\!1$ &$6081$ & $9.3e\!-\!8$ & $1.1e\!-\!7$ & $5.8e\!-\!8$
\\ \hline
 $50$ & $[1.1e\!-\!3, 5.2e\!-\!1]$ & $8.7e\!+\!1$ & $9.5e\!+\!2$& $2621$ & $9.8e\!-\!8$ & $1.1e\!-\!7$ & $1.5e\!-\!7$
\\ \hline
\end{tabular}}
}
\vskip -0.10in
\caption{(Example 1) The integration path $\Gamma$ is the same as in %
Table~\ref{tab:9-3-1}, and $q = 2^4$ in (\ref{equ:11-24-13}). %
For the meanings of the columns, refer to Tables \ref{tab:9-3-1} and %
\ref{tab:9-3-2}.}
\label{tab:9-3-3}
\end{table}


\subsection{Example 2}
\label{Subsec:NumExp:Example2}

\begin{figure}[tbh]
\begin{center}
\includegraphics[width=5.8in]{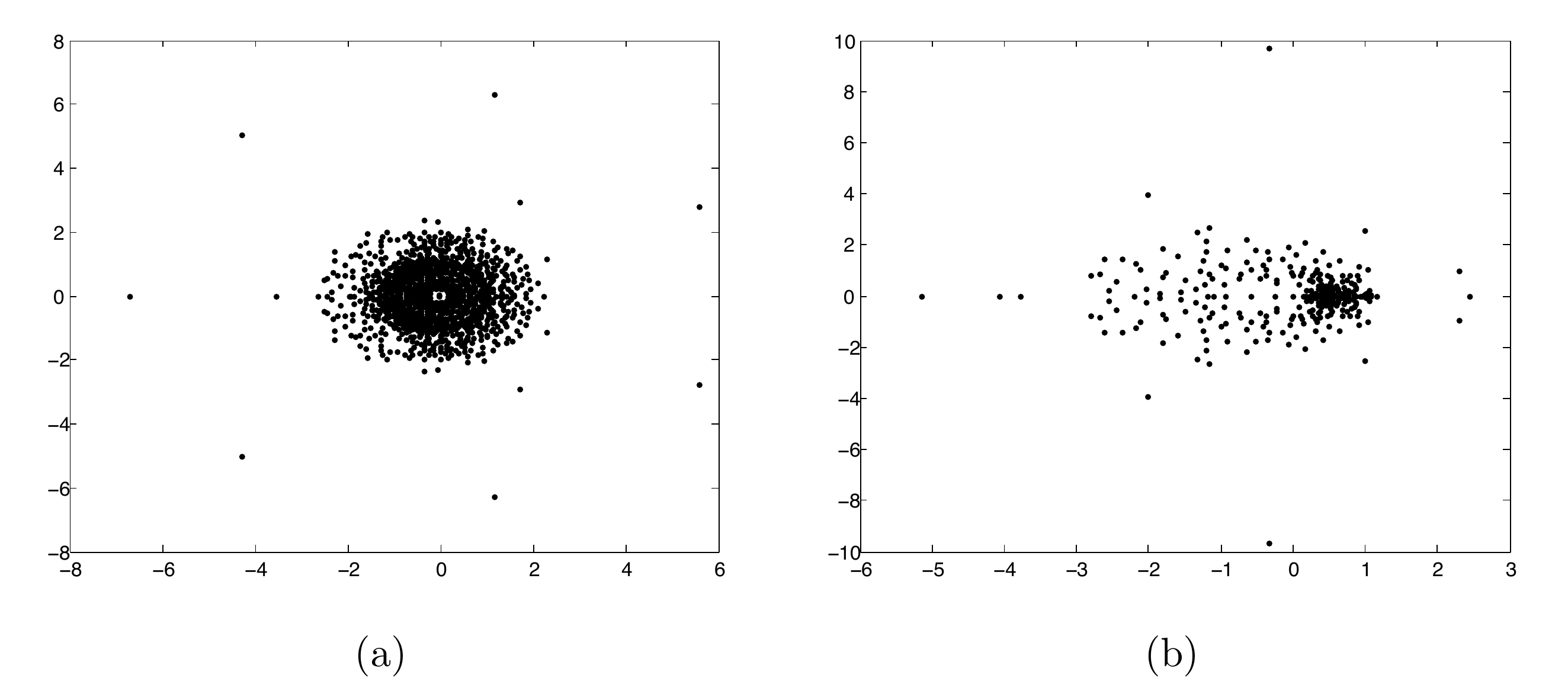}
\end{center}
\vskip -0.25in
\caption{(a) Eigenvalue distribution of the test matrix {\it mahindas}. (b) Eigenvalue distribution of the ILU($0$)-preconditioned {\it mahindas}. }
\label{fig:1-29-1}
\end{figure}

The following two test data sets are part of The University of Florida Sparse
Matrix Collection \cite{ufl_sparse_matrics}.
These data sets have been used in \cite{dly} for the numerical experiments.
\begin{enumerate}
\item[(a)]
{\it bcsstm27} is from a mass matrix buckling problem.
It is a $1224 \times 1224$ real symmetric and indefinite matrix $A$ with
$56126$ nonzero entries.
Per the right hand side in (\ref{equ:11-19-1}), we set $b = A {\bf 1}$,
where ${\bf 1}$ is the vector of all ones.
A spectral plot for {\it bcsstm27} is presented in Figure~\ref{fig:ellipse}(b).
\item[(b)]
{\it mahindas} is from an economics problem.
It is a $1258 \times 1258$ real unsymmetric matrix $A$ with $7682$ nonzero
entries.
Again, we set $b = A {\bf 1}$ as the right hand side in (\ref{equ:11-19-1}).
A spectral plot for {\it mahindas} is in Figure~\ref{fig:1-29-1}(a).
\end{enumerate}

An ILU preconditioner generated by the Matlab function $[L, U, Pr] = luinc(A,
'0')$ was used for {\it mahindas}, namely, instead of solving
(\ref{equ:11-19-1}), we solved

The Matlab function $[L, U, Pr] = luinc(A, '0')$ is used to create an ILU
preconditioner for {\it mahindas}.
Instead of solving (\ref{equ:11-19-1}), we solved
\begin{equation}\label{equ:8-26-1}
\tilde{A} x = \tilde{b}
\end{equation}
by Algorithm \ref{alg:DGMRES}, where $\tilde{A} = L^{-1} Pr A U^{-1}$ and
$\tilde{b} = L^{-1} Pr b$.
The $A$ and $b$ in (\ref{equ:7-31}) are replaced with $\tilde{A}$ and
$\tilde{b}$, respectively.
Once a solution $x$ to (\ref{equ:8-26-1}) is obtained, we compute a solution
to the original system (\ref{equ:11-19-1}) by $U^{-1}x$.
Since the $U$ factor obtained from {\it luinc} has some zeros along its main
diagonal, we replace those zeros by $1$ so that $U$ is invertible.
A spectral plot for $\tilde{A}$ is given in Figure~\ref{fig:1-29-1}(b).

However, we do not apply any preconditioner to {\it bcsstm27}.

In Example~2, we only use MBiCG as the linear solver.
Numerical results are summarized in Tables \ref{tab:11-24-1}-\ref{tab:8-24-2}.

Without deflation MBiCG performs very poorly for both {\it bcsstm27} and the
ILU($0$)-preconditioned {\it mahindas}.
Specifically, MBiCG does not converge within {\it maxit} $= 10^3 N$ iterations
in terms of the relative residues {\it relres2}.
Further, the computed solutions by MBiCG are far from the corresponding exact
solutions $x^* = {\bf 1}$ according to the relative errors {\it relerr}.
With an appropriate eigenvalue-deflation, however, the situation is improved
(see the numerical results in Tables \ref{tab:8-20-1}-\ref{tab:8-24-2}).

The computed matrices $Z$ in Computations \#3 and \#4 in
Table~\ref{tab:8-20-1} worked well for {\it bcsstm27}, but not for the
ILU($0$)-preconditioned {\it mahindas}.
The computed $Z$'s for {\it mahindas} contain some
nearly dependent columns that lead to large condition numbers Cd($M$) of $M$.
The situation is significantly improved after the nearly dependent columns in
$Z$ are removed using Algorithm \ref{alg:8-18-1}.
As a result, MBiCG converged when it solved (\ref{equ:7-31}) (see
Table~\ref{tab:8-24-1}).

Now compare Computations \#4 and \#6 for {\it bcsstm27} in Tables
\ref{tab:8-20-1} and \ref{tab:8-24-1}, respectively.
The two condition numbers Cd($M$) of $M$ do not differ much in magnitude, but
the $m$ ($= 400$) in Computation \#4 is larger than the $m$ ($= 371$) in
Computation \#6.
As a result, MBiCG converged faster in Computation \#4 than it did in
Computation \#6 when solving (\ref{equ:7-31}).

In Computations \#7 and \#8, we randomly picked an initial guess for the
solution of each linear system in (\ref{equ:8-17-1}) and then computed $Z$ by
(\ref{equ:11-24-13}).
The resulting $Z$ has a better performance than the $Z$ obtained with a zero
initial guess as in Computations \#3--\#6.
We first compare Computation \#8 and Computation \#4 for {\it mahindas} in
Tables \ref{tab:8-24-2} and \ref{tab:8-20-1}.
Both condition numbers Cd($M$) of $M$ are about the same in magnitude, but
MBiCG and Algorithm \ref{alg:DGMRES} in Computation \#8 performs much better
than in Computation \#4.

We compare Computations \#8 and \#6 for {\it mahindas} in Tables
\ref{tab:8-24-2} and \ref{tab:8-24-1}.
The column size $m$ of $Z$ in Computation \#8 ($m = 50$) is much larger than
the $m$ in Computation \#6 ($m = 10$).
This explains the faster convergence of MBiCG in Computation \#8 on the
solution of (\ref{equ:7-31}) despite the fact that the $M$ in Computation \#8
is ill conditioned relative to the $M$ in Computation \#6.

Finally, we remark that we have chosen $m \geq s$ in the experiments presented
above.
When $m < s$, Algorithm \ref{alg:DGMRES} plus (\ref{equ:11-24-13}) still works
well, but not as impressively as in the case when $m \geq s$.
For an estimate of $s$, the stochastic method in \cite{futa} should be useful.
See \cite{ycy1} for a concise description of this method.
Moreover, the method in \cite{nps} is also recommended.

The most expensive part in the proposed method of Algorithm \ref{alg:DGMRES}
plus (\ref{equ:11-24-13}) is clearly the computation of $Z$ in
(\ref{equ:11-24-13}).
In \S\ref{sect:future-work}, we describe state of the art parallel multigrid
methods that can be applied to the computation of $Z$.

\begin{table}
\centering
\footnotesize{
{\newcommand{\q}[1]{\mc{1}{|l|}{\small\tt #1}}
\noindent
\begin{tabular}{*{10}{|c}|} \hline
\mc{3}{|c}{}
&
\mc{3}{|c|}{Computation \#1}& \mc{4}{|c|}{Computation \#2}
\\ \hline
& Circle $\Gamma$ & & \mc{3}{|c|}{ Solve (\ref{equ:11-19-1})} &\mc{2}{|c|}{Solve (\ref{equ:7-31})} & \mc{2}{|c|}{Algorithm \ref{alg:DGMRES}} \\ \hline
Matrix & $(c, r)$ & \#eig $\Gamma$
& \#iter & relres2 & relerr & \#iter & relres1 & relres2 & relerr
\\ \hline\hline
 {\it bcsstm27}& $(0, 5)$ &$363$
   &$1224000$  &$2.0 e\!-\!6$ & $5.7e\!-\!1$ & $227$ & $9.5e\!-\!8$& $9.5e\!-\!8$& $2.4e\!-\!7$
\\ \hline
{\it mahindas}&$(-1,1)$ &$31$ &
   $1258000$  &$3.5e\!-\!7$ & $1.8$ & $1841$ &$5.6e\!-\!8$ &$5.6e\!-\!8$& $4.4e\!-\!3$
\\ \hline
\end{tabular}}
}
\vskip -0.10in
\caption{(Example 2)
The linear systems
(\ref{equ:11-19-1}) and (\ref{equ:7-31}) were solved by MBiCG.
$\Gamma$ is a circle with center $c$ and radius $r$.
{\it \#eig $\Gamma$} stands for the number of eigenvalues of $A$ inside $\Gamma$,
and {\it \#iter} the number of iterations by MBiCG.
For {\it mahindas}, an ILU($0$) preconditioner was applied.
}
\label{tab:11-24-1}
\end{table}

\begin{table}
\centering
\footnotesize{
{\newcommand{\q}[1]{\mc{1}{|l|}{\small\tt #1}}
\noindent
\begin{tabular}{*{9}{|c}|} \hline
\mc{9}{|c|}{Computation \#3}
\\ \hline
& \mc{2}{|c|}{Solve (\ref{equ:8-17-1})} & \mc{2}{|c|}{Compute $Z$ by (\ref{equ:11-24-13})} &\mc{2}{|c|}{Solve (\ref{equ:7-31})} & \mc{2}{|c|}{Algorithm \ref{alg:DGMRES}}
\\ \hline
Matrix &$m$& Rel. Res. Range & Cd($Z$) & Cd($M$)& \#iter & relres1 & relres2 & relerr
\\ \hline 
 {\it bcsstm27}& $400$& $[4.1e\!-\!10, 6.1e\!-\!2]$ & $1.5e\!+\!3$&$1.0e\!+\!6$ &$843$ &$6.3e\!-\!8$ &$6.2e\!-\!8$ &$5.0e\!-\!6$
\\ \hline
{\it mahindas}&$50$&$[1.2e\!-\!10, 1.0]$ &$1.3e\!+\!2$ &$1.7e\!+\!6$  &$12580$ & $6.4e\!-\!2$ & $2.2e\!-\!1$ &$2.5e\!+\!4$\\ \hline
 \mc{9}{|c|}{Computation \#4}
\\ \hline
 {\it bcsstm27}&$400$& $[3.2e\!-\!14, 8.0e\!-\!4]$ &$2.5e\!+\!3$ &$4.5e\!+\!6$ &$563$ &$9.2e\!-\!8$ &$9.2e\!-\!8$ &$3.4e\!-\!2$
\\ \hline
{\it mahindas}&$50$& $[1.5e\!-\!10, 1.0]$  & $2.5e\!+\!2$ &$2.2e\!+\!6$ &$12580$ & $3.1e\!-\!2$ & $3.6e\!-\!2$ &$1.3e\!+\!4$
\\ \hline
\end{tabular}}
}
\vskip -0.10in
\caption{(Example 2) The integration paths $\Gamma$ are shown in Table~\ref{tab:11-24-1}, and $q = 2^4$ in (\ref{equ:11-24-13}).
The linear systems
(\ref{equ:8-17-1}) and (\ref{equ:7-31}) were solved by MBiCG.
For the meanings of the columns, refer to Tables \ref{tab:9-3-1} and \ref{tab:9-3-2}.
}
\label{tab:8-20-1}
\end{table}

\begin{table}
\centering
\footnotesize{
{\newcommand{\q}[1]{\mc{1}{|l|}{\small\tt #1}}
\noindent
\begin{tabular}{*{10}{|c}|} \hline
\mc{10}{|c|}{Computation \#5}
\\ \hline
& \mc{2}{|c|}{Solve (\ref{equ:8-17-1})} & \mc{3}{|c|}{Algorithm \ref{alg:8-18-1}} &\mc{2}{|c|}{Solve (\ref{equ:7-31})} & \mc{2}{|c|}{Algorithm \ref{alg:DGMRES}}
\\ \hline
Matrix &$m$ & Rel. Res. Range & {\it rk} of $Z$& Cd($Z$) & Cd($M$) & \#iter & relres1 & relres2 & relerr
\\ \hline 
 {\it bcsstm27}& $400$& $[4.1e\!-\!10, 6.1e\!-\!2]$ &$370$ &$9.3e\!+\!1$ &$1.2e\!+\!5$ &$2824$ &$8.8e\!-\!8$ &$8.8e\!-\!8$ &
$6.2e\!-\!2$\\ \hline
{\it mahindas}& $50$ & $[1.2e\!-\!10, 1.0]$  &$12$ &$2.7e\!+\!1$ &$8.9e\!+\!2$ & $9318$ & $8.6e\!-\!8$ &$8.6e\!-\!8$ &$6.0e\!-\!3$
\\ \hline \hline
 \mc{10}{|c|}{Computation \#6}
\\ \hline
 {\it bcsstm27}&$400$ & $[3.2e\!-\!14, 8.0e\!-\!4]$ &$371$ &$9.4e\!+\!1$ &$1.9e\!+\!5$ &$2224$ &$5.8e\!-\!8$ &$5.8e\!-\!8$ &
$6.3e\!-\!2$\\ \hline
{\it mahindas}&$50$ & $[1.5e\!-\!10, 1.0]$  &$10$ &$1.9e\!+\!1$ &$4.4e\!+\!2$ & $4688$ & $6.8e\!-\!8$ &$6.8e\!-\!8$ &
$3.5e\!-\!3$
\\ \hline
\end{tabular}}
}
\vskip -0.10in
\caption{(Example 2) The integration paths $\Gamma$ are shown in Table~\ref{tab:11-24-1}, and $q = 2^4$ in (\ref{equ:11-24-13}). The linear systems
(\ref{equ:8-17-1}) and (\ref{equ:7-31}) were solved by MBiCG.
The numbers in the column titled ``{\it rk} of $Z$'' are the
numerical ranks {\it rk} of $Z$ output by Algorithm \ref{alg:8-18-1}. For the meanings of other columns,
refer to
Tables \ref{tab:9-3-1} and \ref{tab:9-3-2}.
}
\label{tab:8-24-1}
\end{table}

\begin{table}
\centering
\footnotesize{
{\newcommand{\q}[1]{\mc{1}{|l|}{\small\tt #1}}
\noindent
\begin{tabular}{*{9}{|c}|} \hline
\mc{9}{|c|}{Computation \#7}
\\ \hline
& \mc{2}{|c|}{Solve (\ref{equ:8-17-1})} & \mc{2}{|c|}{Compute $Z$ by (\ref{equ:11-24-13})} &\mc{2}{|c|}{Solve (\ref{equ:7-31})} & \mc{2}{|c|}{Algorithm \ref{alg:DGMRES}}
\\ \hline
Matrix &$m$&  Rel. Res. Range & Cd($Z$) & Cd($M$)& \#iter & relres1 & relres2 & relerr
\\ \hline
 {\it bcsstm27}& $400$& $[3.9e\!-\!9, 2.2e\!-\!1]$ &$5.8e\!+\!2$ &$5.4e\!+\!5$ &$1114$ &$8.2e\!-\!8$ &$8.4e\!-\!8$ &$1.9e\!-\!6$
\\ \hline
{\it mahindas}& $50$& $[2.8e\!-\!10, 2.0e\!+1]$&$8.5e\!+\!1$ & $1.2e\!+\!5$ &$12580$ & $5.6e\!-\!5$ & $5.7e\!-\!5$ &$4.8$
\\ \hline \hline
 \mc{9}{|c|}{Computation \#8}
\\ \hline
 {\it bcsstm27}&$400$ & $[2.7e\!-\!13, 3.1e\!-\!3]$& $2.3e\!+\!3$ &$4.3e\!+\!6$ &$686$ &$9.1e\!-\!8$ &$9.1e\!-\!8$ &$3.2e\!-\!6$
\\ \hline
{\it mahindas}& $50$&$[9.8e\!-\!11, 1.3e\!+\!1]$&$6.7e\!+\!2$ &$9.7e\!+\!6$  &$2286$ & $8.4e\!-\!8$ & $8.5e\!-\!8$ &$1.8e\!-\!2$\\ \hline
\end{tabular}}
}
\vskip -0.10in
\caption{(Example 2) The integration paths $\Gamma$ are shown in Table~\ref{tab:11-24-1}, and $q = 2^4$ in (\ref{equ:11-24-13}).
The linear systems
(\ref{equ:8-17-1}) and (\ref{equ:7-31}) were solved by MBiCG.
For the meanings of the columns,
refer to
Tables \ref{tab:9-3-1} and \ref{tab:9-3-2}.
}
\label{tab:8-24-2}
\end{table}



\section{Future Work}
\label{sect:future-work}

We can formulate either geometric multigrid
\cite{GPAstrakhantsev_1971a,
NSBakhvalov_1966a,
ABrandt_1977b,
CCDouglas_1984a,
RPFedorenko_1961a,
RPFedorenko_1964a,
WHackbusch_1985a,
PWesseling_1992a}
or algebraic multigrid 
\cite{KStuben_2000a}
using the same
notation level to level using the abstract multigrid approach developed in
\cite{CCDouglas_1984c,
CCDouglas_1995a,
CCDouglas_JDouglas_DEFyfe_1994a,
REBank_CCDouglas_1985a,
CCDouglas_1994b,
CCDouglas_1995a,
CCDouglas_JDouglas_DEFyfe_1994a}.

Assuming the cost of the smoother (or rougher) on each level is $O(N_j)$,
$j=1,\cdots,k$, Algorithm MGC with $p$ recursions to solve problems on level
$k-1$ has complexity
$$
\label{eqn:multigrid-complexity}
W_{MGC}(N_k) =
\left\{
\begin{array}{ll}
O(N_k)          & 1 \leq p \leq \sigma\\
O(N_k \log N_k) & p = \sigma\\
O(N_k^{\log p}) & p > \sigma.\\
\end{array}
\right.
$$
Under the right circumstances, multigrid is of optimal order as a solver.

Consider the example (\ref{equ:11-24-10}) in \S\ref{Subsec:Problem:NumerExI}.
A simple geometric multigrid approximation to (\ref{equ:11-24-10}) produces
a very good solution in $4$~V Cycles or $2$~W cycles using
the deflated GMRES as the rougher.
Each V or W Cycle is $O(N_k)$.
Hence, we have an optimal order solver for (\ref{equ:11-24-10}), which
would not be the case if we used BiCG or deflated GMRES on a single grid.

High performance computing versions of multigrid based on using
hardware acceleration with memory caches was extensively studied in the
early 2000's
\cite{CCDouglas_JHu_MKowarschik_URude_CWeiss_2000a}.

Parallelization of Algorithm MGC is straightforward
\cite{CCDouglas_1996c}.
\begin{itemize}
\item
For geometric multigrid, on each level~$j$,
data is split using a
domain decomposition paradigm.
Parallel smoothers (roughers) are used.
The convergence rate degrades from the standard serial theoretical rate,
but not by a lot, and scaling is good given sufficient data.
\item
For algebraic multigrid, the algorithms can be either straightforward
(e.g., Ruge-Studen \cite{JWRuge_KStuben_1985a}
or
Beck \cite{RBeck_1999a}) to quite complicated (e.g.,
AMGe \cite{GHaase_2000a}).
Solutions have existed for a number of years, so it is a matter of
choosing an exisiting implementation.
In some cases, using a tool like METIS or ParMETIS is sufficient to
create a domain decomposition-like system based on graph connections in
$A_j$, which reduces parallelization back to something similar to the
geometric case.
\end{itemize}
In many cases, the complexity of this type of parallel multigrid for $P$
processors becomes
$$
\label{eqn:ParallelMGCcost}
W_{MGC,P}(N_k) = W_{MGC}(N_k)\log P / P.
$$


\section{Conclusions}
\label{sect:conclusions}

We incorporate the delation projector $P$ in (\ref{equ:8-21-1}), with $Z$
defined by and computed by (\ref{equ:11-24-1}) and (\ref{equ:11-24-13}),
respectively, into Krylov subspace methods to enhance the stability and
accelerate the convergence of the iterative methods for solving ill
conditioned linear systems.
Our experiments suggest that Algorithm~\ref{alg:DGMRES} plus
(\ref{equ:11-24-13}) has the potential to solve ill conditioned problems much
faster and more accurately than standard Krylov subspace solvers.
Moreover, to our best knowledge, the constructions of most, if not all,
deflation subspace matrices $Z$ in the literature are problem dependent.
The method proposed here, however, is problem independent.

More experiments, especially on test data of large size (e.g., millions or
more unknowns) are needed to better understand the behavior of the proposed
algorithm.
Implementation of robust and efficient parallel multigrid methods for solving
(\ref{equ:8-17-1}) and the realization of a software package for a wide
variety of applications is currently under our investigation.


\section{Appendix}\label{sect:appendix}

The following Algorithm \ref{alg:8-18-1} employs the Gaussian elimination with
complete pivoting (CGE) to detect the rank {\it rk} and to select linearly
independent columns of an input $Z \in {\mathbb C}^{N \times m}$, where
$\alpha$ is a comparison parameter and {\it tol\_cge} is a stopping tolerance.
The output $Z$ of the algorithm is a $N$-by-$rk$ matrix consisting of the selected columns of the input $Z$.

\begin{algorithm}\label{alg:8-18-1}
  \SetKwData{Left}{left}\SetKwData{This}{this}\SetKwData{Up}{up}
  \SetKwFunction{Union}{Union}\SetKwFunction{FindCompress}{FindCompress}
  \SetKwInOut{Input}{input}\SetKwInOut{Output}{output}
\begin{tabbing}
x\=xxx\= xxx\=xxx\=xxx\=xxx\=xxx\kill
\>Function $[Z, rk] = \mbox{cge}(Z, \alpha, tol\_cge)$ \\
\>1. \> Form $\widehat{Z} = Z^H Z \in {\mathbb C}^{m \times m}$; Set $rk = m$.\\
\>2. \> Determine $(i_0, j_0)$ with $1 \leq i_0, j_0 \leq m$ so that\\
\>\>\>\> $|\widehat{Z}_{i_0 j_0}| = \max\{|\widehat{Z}_{ij}|: 1 \leq i, j \leq m\}$.\\
\>3.\> If $|\widehat{Z}_{i_0 j_0}| < \alpha$, then set $rk = 0$ and $Z = [\,\,]$; Stop.\\
\>4.\> Set $\alpha = |\widehat{Z}_{i_0 j_0}|$;\\
\>5.\> $\widehat{Z}_{(:,1)} \leftrightarrow \widehat{Z}_{(:,j_0)}$;
$\widehat{Z}_{(1,:)} \leftrightarrow \widehat{Z}_{(i_0,:)}$. \,\, \% move $\widehat{Z}_{i_0j_0}$ to the $(1,1)$ position.\\
\>6.\>$Z_{(:,1)} \leftrightarrow Z_{(:,j_0)}$.\,\, \% interchange the $1$st and the $j_0$th columns of $Z$.\\
\>7.\> For $j = 1:m - 1$\\
\>8.\>\>For $i = j+1:m$\\
\>9.\>\>\> $\widehat{Z}_{(i,j:m)} = \widehat{Z}_{(i,j:m)} -
(\widehat{Z}_{ij}/\widehat{Z}_{jj})
\widehat{Z}_{(j,j:m)}$. \,\, \% perform elimination.\\
\>10.\>\>End\\
\>11.\>\>Determine $(i_0, j_0)$ with $j+1 \leq i_0, j_0 \leq m$ so that\\
\>\>\>\>\> $|\widehat{Z}_{i_0 j_0}| = \max\{|\widehat{Z}_{p q}|: j+1 \leq p,
q \leq m\}$.  \\
\>12.\>\> If $|\widehat{Z}_{i_0 j_0}|/\alpha < tol\_cge$, set $rk = j$ and $Z =
Z_{(:,1:rk)}$; Stop.\\
\>13.\>\>$\widehat{Z}_{(:,j+1)} \leftrightarrow \widehat{Z}_{(:,j_0)}$;
$\widehat{Z}_{(j+1,:)} \leftrightarrow \widehat{Z}_{(i_0,:)}$.\,\, \% move $\widehat{Z}_{i_0j_0}$ to the $(j+1, j+1)$ position.\\
\>14.\>\>$Z_{(:,j+1)} \leftrightarrow Z_{(:,j_0)}$.\,\, \% interchange column $j_0$ and column $j+1$ of $Z$.\\
\>15.\> End
\end{tabbing}
\caption{Gaussian elimination with complete pivoting to detect the rank and to select linearly independent columns of an input matrix $Z$.}
\end{algorithm}

In Line 3 of Algorithm \ref{alg:8-18-1},
we consider $\widehat{Z} = 0$ when $|\widehat{Z}_{i_0 j_0}|$ is small compared
to $\alpha$, and hence the rank {\it rk} of $Z$ is zero.
Similarly, in Line 12, if $|\widehat{Z}_{i_0 j_0}|$ is small compared to
$\alpha$, then we regard $\widehat{Z}_{(j+1:m,j+1:m)}$ as $0$, and therefore
$rk = j$.


\section*{Acknowledgments}
\label{sect:acks}

This research was supported in part by National Science Foundation grants
ACI-1440610, ACI-1541392, and DMS-1413273.
In addition, we would like to thank
Iain Duff,
Ken Hayami,
and
Dongwoo Sheen
for their valuable comments and suggestions for future work.


%
\label{sect:bib}
\bibliographystyle{plain}
\bibliography{submit}







\end{document}